\documentclass[leqno]{article}
\usepackage{amsmath, amssymb}
\def\beq{\begin{equation}}
\def\eeq{\end{equation}}
\def\bea{\begin{eqnarray*}}
\def\eea{\end{eqnarray*}}
\def\b{\begin{eqnarray}}
\def\e{\end{eqnarray}}
\def\a{\begin{align}}
\def\enda{\end{align}}

\numberwithin{equation}{section}
\newtheorem{theorem}{T\begin{scriptsize}HEOREM\end{scriptsize}}[section]
\newtheorem{corollary}[theorem]{C\begin{scriptsize}OROLLARY\end{scriptsize}}
\newtheorem{lemma}[theorem]{L\begin{scriptsize}EMMA\end{scriptsize}}
\newtheorem{proposition}[theorem]{P\begin{scriptsize}ROPOSITION\end{scriptsize}}

\def\na{\nabla}
\def\pa{\partial}

\def\pt{\pa_\theta}
\def\ptt{\langle \pt \rangle}

\def\operatorname{\mbox}

\def\Tr{\operatorname{Tr}}
\def\div{\operatorname{div}\,}
\def\de{\operatorname{def}\,}
\def\curl{\operatorname{curl}\,}

\def\R{{\bf R}}
\def\Om{\Omega}
\def\pOm{\pa \Om}
\def\ve{\varepsilon}
\def\bp{{\bf \noindent P\begin{scriptsize}ROOF\end{scriptsize}: }}
\newcommand{\ep}{\hfill\vrule height6pt width6pt depth0pt}

\begin{document}

\title{{\it A priori} estimates for the motion of a self-gravitating incompressible liquid with
free surface boundary.}
\author{Hans Lindblad and Karl H\aa kan Nordgren}
\date{\today}
\maketitle

\begin{abstract} In this paper, we prove {\it a priori} estimates in Lagrangian coordinates for the equations of motion of an incompressible, inviscid, self-gravitating fluid with free boundary. The estimates show that on a finite time interval we control five derivatives of the fluid velocity and five and a half derivatives of the coordinates of the moving domain.
\end{abstract}

\section{Introduction.}

Let $\Om_t \subseteq \R^n$ be the domain occupied by a fluid at time
$t \in [0,T]$ and suppose that the fluid has velocity $v(t,x)$ and
pressure $p(t,x)$ at a point $x$ in $\Om_t$. For an inviscid,
self-gravitating fluid these two quantities are related by Euler's
equation \beq \label{eq:DebussyTroisPoemesMov1} \big(\pa_t +
v^i\pa_i\big)v_j = - \pa_j p - \pa_j \phi \eeq in $\Om_t$, where
$\pa_i = \frac{\pa}{\pa x^i}$ and $v^i = \delta^{ij}v_j$ and where
$\phi$ is the Newtonian gravity-potential defined by \beq
\label{eq:mar25080754} \phi(t,x) = - \chi_{\Om_t} * \Phi (x)\eeq on
$\Om_t$, where $\chi_{\Om_t}$ is a function which takes the value
$1$ on $\Om_t$ and the value $0$ on the complement of $\Om_t$ and
where $\Phi$ is the fundamental solution to the Laplacean. Thus
$\phi$ satisfies $\Delta \phi = - 1$ on $\Om_t$. We can impose the
condition that the fluid be incompressible by requiring that the
fluid-velocity be divergence-free: \beq
\label{eq:BarberViolinConcerto1Mov2} \div v = \pa_i v^i = 0 \mbox{
in } \Om_t. \eeq The absence of surface-tension is imposed with the
following boundary condition: \beq
\label{eq:BarberViolinConcerto1Mov2b} p = 0 \mbox{ on } \pa \Om_t
\eeq and to make the free-boundary move with the fluid-velocity, we
have \beq \label{eq:CesarFranckPreludeAriaEtFinal} \pa_t + v^i\pa_i
\mbox{ is in the tangent-space of } \cup_{t \in [0,T]} [\Om_t\times
\{t\}]. \eeq The problem is, then, to prove {\it a priori} estimates
for $v$ satisfying (\ref{eq:DebussyTroisPoemesMov1}) -
(\ref{eq:CesarFranckPreludeAriaEtFinal}) in some interval $[0,T]$,
given the initial-conditions \beq
\label{eq:BarberViolinConcerto1Mov3b} v = v_0 \mbox{ on } \Om_0,
\eeq where $v_0$ and $\Om_0$ are known. We will also assume
that initially there is a constant $c_0$ such that \beq
\label{eq:BarberViolinConcerto1Mov3}  \na p \cdot N \leq -c_0 < 0
\mbox{ on } \pOm_0 \eeq where $N$ is the exterior unit normal to
$\pOm$, \eqref{eq:BarberViolinConcerto1Mov3} is a natural physical
condition since the pressure of a fluid has to be positive and the
problem is ill-posed if this is not satisfied, see
Ebin \cite{MR886344}. This condition is related to Rayleigh-Taylor
instability.
\bigskip

We will assume for simplicity that $n$, the number of
space-dimensions, is $2$. We will also assume that there is a volume-preserving diffeomorphism $x_0: \Om \to \Om_0$, where $\Om = \{y \in \R^2: |y| < 1\}$, which will allow us describe the derivatives which are
tangential to $\Om_t$ in a particularly simple way, and it
will also allow fractional derivatives to be defined by using
Fourier series without recourse to partitions of unity. That part of
the argument can be used, with minor modifications, in the case of
arbitrary space dimension. The dimension will also allow simpler
energy estimates because the curl of the velocity is conserved.

Suppose now that $v$ satisfies (\ref{eq:BarberViolinConcerto1Mov2}). We define Lagrangian coordinates as follows: Define $x$ by \beq
\label{eq:mar25081050} \frac{dx}{dt}(t,y) = v(t,x(t,y)) \mbox{ and }
x(0,y) = x_0(y) \eeq for $y$ in $\Om$ and for $t$ in some time interval
$[0,T]$. Since $v$ satisfies
(\ref{eq:BarberViolinConcerto1Mov2}) and this means that \beq \label{eq:mar25080910} \pa_t
\det\left(\frac{\pa x}{\pa y}\right)(t,y) = \div v \circ x(t,y) = 0.
\eeq And since $\det\left(\frac{\pa x}{\pa y}\right)(0,y) = 1$ we
therefore have $\det\left(\frac{\pa x}{\pa y}\right) = 1$ in $\Om$.
We will prove the following theorem:

\begin{theorem} \label{may16080841} Let $v$ satisfy (\ref{eq:DebussyTroisPoemesMov1})
and (\ref{eq:BarberViolinConcerto1Mov2}) and let $p$ satisfy
(\ref{eq:BarberViolinConcerto1Mov2b}) and
(\ref{eq:BarberViolinConcerto1Mov3}). Let the flow $x$ of $v$ be
defined by (\ref{eq:mar25081050}) and define $V(t,y) = v(t,x(t,y))$.
Define \[ E(t) = \sup_{[0,t]} \left[\|V\|_{5} + \|x\|_{5.5} +
\|\curl(v)\|_{H^{4.5}(\Om_t)}\right]. \] Then there is $T > 0$ such that $E(T) \leq P\big[E(0)\big]$ where $P$ is a
polynomial. \end{theorem}

\subsection{Background.}

Past progress has been made in three situations: The first progress was made on the water-wave
problem under the assumption that the fluid be irrotational $-$ that is, the curl of the
fluid-velocity is zero $-$, incompressible and that the free-boundary not
be subject to surface-tension. Notable results in this area are Wu's papers
 \cite{MR1471885} and \cite{MR1641609} where she uses Clifford analysis to show
 well-posedness in two and then three dimensions in an infinitely deep fluid; and also
 Lannes' paper \cite{MR2138139} where the Nash-Moser technique is used to prove
 well-posedness in arbitrary space-dimesions for a fluid of finite depth.

In \cite{MR1780703}, Christodoulou and Lindblad proved {\it a priori} estimates for
the incompressible Euler's equation, without the assumption of irrotationality.
They were not sufficient to obtain the existence result, however, because no
 approximation-schemes was discovered which did not destroy the structure in the
 equations on which the estimates relied. In \cite{MR1934619} Lindblad proved that the linearized
 equations are well-posed.
 In \cite{MR2178961} Lindblad then used
 the Nash-Moser approximation scheme to obtain the full
 well-posedness.
  Well-posedness was also proved by Coutand and Shkoller in
 \cite{MR2291920}, using a fixed-point argument which relies on smoothing the
 fluid-velocity only $-$ crucially $-$ in the direction tangential to the boundary.
 This is followed by energy estimates which we will discuss in detail in section
 \ref{mar26080821}. Also, in \cite{SZ1}, Shatah and Zeng prove {\it a priori}
 estimates under these conditions by considering Euler's equation as the geodesic
 equation on the group of volume-preserving diffeomorphisms. The latter two papers
 also consider the case of positive surface-tension.

\section{Preliminaries.}

We will let $x$ be coordinates on $\Om_t$ defined by
(\ref{eq:mar25081050}) and let $\pa_i, \pa_j, \pa_k, \ldots$ be
derivatives on $\Om_t$; and we will let $y$ be coordinates on $\Om$
and let $\pa_a, \pa_b, \pa_c, \ldots$ be derivatives on $\Om$. Also,
we will let $\na$ denote an arbitrary derivative on $\Om_t$ and
$\pa$ be an arbitrary derivative on $\Om$.

\subsection{Change of variables.}

Let $A^a_i(t,y) = \frac{\pa y^a}{\pa x_i}\circ x(t,y)$ and let $B_a^i(t,y) = \frac{\pa x_i}{\pa y^a}(t,y)$.
By (\ref{eq:mar25080910}) we see that $\det(B) = 1$ and hence $\det(A) = 1$ as well in a
time interval $[0,T]$. This means that $x(t,\cdot):\Om \to \Om_t$ is a change of variables.

\subsection{Norms.} \label{mar25081312}

Let $f:\Om \to \R^2$. Define $\|f\|^2 = \int_{\Om} \delta_{ij} f^i(y)f^j(y) dy$ and for
an integer $k \geq 0$ we define \[\|f\|_{k}^2 = \sum_{i = 0}^k \|\na^i[f]\|^2 \] where
$\na = \left(\frac{\pa}{\pa y^1}, \frac{\pa}{\pa y^2}\right)$. We define the
intermediate spaces by interpolation, see for instance \cite{MR0350177} and
\cite{MR0450957}. For a function $g: \Om_t \to \R^2$, define
$\|g\|_{L^2(\Om_t)}^2 = \int_{\Om_t} \delta_{ij} g^i(x)g^j(x) dx$
and for an integer $k \geq 0$ define
\[\|g\|_{H^k(\Om_t)}^2 = \sum_{i = 0}^k \|\na^i[g]\|_{L^2(\Om_t)}^2 \] where
$\na = \left(\frac{\pa}{\pa x^1}, \frac{\pa}{\pa x^2}\right)$.
Again, we define the intermediate spaces by interpolation.
For a function $h: \pOm_t \to \R^2$, define
$\|h\|_{L^2(\pOm_t)}^2 = \int_{\pOm_t} \delta_{ij} h^i(x)h^j(x) dS(x)$

\subsection{Special derivatives.} \label{mar25080913}

In this paper, we will make use of two special derivatives. First we have
derivatives which are tangential to the boundary:
\[ y^a \pa_b - y^b \pa_a \mbox{ for } a,b = 1,2 \]
where the summation convention is {\it not} employed.
We will abuse notation and denote these derivatives and their
push-forwards on $\Om_t$ with $\pt$. The second type of special
derivative is defined as follows: Let $f: \Om \to \R^2$ and let's
abuse notation and write $f(\rho,\theta)$ to mean the polar representation of $f$.
Let $\sum f_k(\rho) e^{ik \theta}$ be the tangential Fourier expansion of $f(\rho,\theta)$.
We now define a tangential Sobolev-type-derivative $\ptt^s$ to be an operator which sends
\[\sum f_k(\rho) e^{ik \theta} \mbox{ to } \sum \langle k \rangle^s f_k(\rho) e^{ik \theta}, \]
where $\langle k \rangle  = \left[1 + |k|^2\right]^\frac12$.
For a function $g$ on $\Om_t$, define
$\|\ptt^s[g]\|_{L^2({\Om_t})} = \|\ptt^s[g \circ x ]\|_{L^2(\Om)}$.
For integral $s$ the operator $\ptt^s$ is equivalent
(in the $L^2(\Om)$- and $L^2(\pOm)$-norm) to the application of a collection of $\pt$.
Finally, for a function $h: \pOm_t \to \R^2$, define
 \[\|h\|_{H^s(\pOm_t)}^2 = \|\ptt^s h\|_{L^2(\pOm_t)}^2 \] for a real number $s$.

\subsection{Cut-off functions.} \label{april02081714}

Fix $d_0$ such that the normal $N$ to $\pOm_t$ can be extended into
the image of the set $\{y \in \R^2: 1-d_0 < |y|\}$ under $x$. This
fact is used in lemma \ref{normalderivnoninv} which is presented in
section \ref{april29081259}. Let $\eta_i$ and $\zeta_i$ be radial
functions which form a partition of unity subordinate to the sets
$\{y \in \R^2: |y| < 1-\frac{d_0}{2i}\}$ and $\{y \in \R^2:
1-\frac{d_0}{i} < |y|\}$ respectively. This means that $\eta_i$
takes the value $1$ on the set $\{y \in \R^2: |y| \leq
1-\frac{d_0}{i}\}$ and $\zeta_i$ takes the value $1$ on the set $\{y
\in \R^2: 1-\frac{d_0}{2i} \leq |y|\}$. We will also let $\eta_i$
and $\zeta_i$ denote the analogous functions in the Eulerian frame.

\subsection{Hodge-decomposition inequalities.}\label{april29081259}

In this section we present two divergence-curl estimates which are
used throughout this text. The first allows point-wise control on
all derivatives near the boundary of $\Om_t$ by the divergence, the
curl and tangential derivatives. Letting $\zeta = \zeta_i$ we have
the following:

\begin{lemma} \label{normalderivnoninv} Let $\alpha$ be a vector-field on $\Om_t$.
 Define $(\curl \alpha)_{jk}=\pa_j \alpha_k-\pa_k \alpha_j$ and $\div \alpha = \pa_j \alpha^j$.
 Then we have the following point-wise estimate on $\Om_t$:
 \beq |\zeta \na \alpha| \leq |\zeta \curl \alpha| + |\zeta \div \alpha| + |\zeta \,\pt \alpha|,\eeq
 where $|\cdot|$ denotes the usual Euclidean distance. \end{lemma}

Using lemma \ref{normalderivnoninv} and an induction argument we have the following lemma:

\begin{lemma} \label{april04080805} For $1 \leq s \leq 5$,
\beq \|\zeta \alpha\|_{H^s(\Om_t)} \leq P\big[\|x\|_5\big]
\left[\|\zeta \alpha\|_{L^2(\Om_t)} +  \|\zeta \curl \alpha\|_{H^{s-1}(\Om_t)}
+ \|\zeta \div \alpha\|_{H^{s-1}(\Om_t)} + \sum_{j=1}^s\|\zeta \pt^j
\alpha\|_{L^2(\Om_t)}\right]. \eeq \end{lemma}

We will also use the following estimates which allows $H^s(\Om_t)$ control in terms of the divergence,
 the curl and boundary derivatives:

\begin{lemma} \label{CSineq} Let $\div \alpha$ and $\curl \alpha$ be defined as in lemma
\ref{normalderivnoninvII}. Then, for $s \leq 5$,
\begin{align} \|\alpha\|_{H^s(\Om_t)} \leq P\big[\|x\|_5\big]\left[\|\alpha\|_{L^2(\Om_t)}
+ \|\div \alpha\|_{H^{s-1}({{\Om_t}})} + \|\curl \alpha\|_{H^{s-1}({{\Om_t}})}
+ \|(\ptt^{s-\frac12} \alpha) \cdot N\|_{L^2(\pOm_t)}\right], \end{align}
where $N$ is the outward unit normal to $\pOm_t$. Also, for $s \leq 5$,
\begin{align} \|\alpha\|_{H^s(\Om_t)} \leq  P\big[\|x\|_5\big]
\left[ \|\alpha\|_{L^2({{\Om_t}})} + \|\div \alpha\|_{H^{s-1}(\Om_t)}
+ \|\curl \alpha\|_{H^{s-1}(\Om_t)} + \|(\ptt^{s-\frac12} \alpha)
\cdot Q\|_{L^2(\pOm_t)}\right] \end{align} where $Q$ is a unit vector which is
tangent to $\pOm_t$. \end{lemma}

\section{Elliptic estimates.}

In section \ref{mar26080821}, we will prove energy estimates for
(\ref{eq:DebussyTroisPoemesMov1}) - (\ref{eq:BarberViolinConcerto1Mov3}) and to prepare for this,
we need the elliptic estimates for $p$ and $\phi$ contained in this section.

\subsection{Estimates for $\phi$.}

In this section we prove the following theorem:

\begin{theorem} \label{ChopinEtudeAFlat}
$\|\na \phi\|_{H^5(\Om_t)} \leq P\big[\|x\|_{5.5}\big]$,
where $P$ is a polynomial. \end{theorem}

Using the cut-off functions defined in section \ref{april02081714},
we have $\|\na \phi\|_{H^5(\Om_t)} \leq \|\eta_1 \na \phi\|_{H^5(\Om_t)}
+ \|\zeta_1 \na \phi\|_{H^5(\Om_t)}$.
 This allows us to consider interior and boundary regularity separately.

\subsubsection{Interior regularity.}

In this section we prove the following result:

\begin{theorem} \label{april29081232} For any $1 \leq i$ and $1 \leq s \leq 5$
\begin{align} \label{eq:april29081315} \|\na^s [\eta_i \na \phi]\|_{L^2(\Om_t)}
\leq P\left[\|x\|_{5}\right] \end{align} where $P$ is a polynomial.\end{theorem}

We prove that (\ref{eq:april29081315}) holds by induction on $s$. Suppose that $s=0$.
 We have
 $\|\eta_i\na\phi\|_{L^2(\Om_t)} \leq \|\eta_i\|_{L^\infty(\Om_t)} \\ \times \|\na\phi\|_{L^2(\Om_t)}$
 and \begin{align} \|\na\phi\|_{L^2(\Om_t)}^2 \label{eq:april29081348}
 = \int_{\Om_t} (\pa_j\phi)(\pa^j\phi)dx = \int_{\pOm_t} N^j(\pa_j\phi)\,\phi\,dx
 - \int_{\Om_t} \Delta \phi\,\phi\,dx. \end{align} Since we have
 $\|\phi\|_{L^\infty(\Om_t)} \leq P\big[\|x_\kappa\|_2\big]$ and
 $\|\na \phi\|_{L^\infty(\Om_t)} \leq P\big[\|x_\kappa\|_2\big]$, we control both terms
 in (\ref{eq:april29081348}) appropriately. This proves (\ref{eq:april29081315})
  for $s=1$. Now suppose that $s = 5$ and that we have the result for smaller $s$.
  Then \begin{align} \|\na^5 [\eta_i \na \phi]\|_{L^2(\Om_t)}^2 \label{eq:may14080937}
   & = \int_{\Om_t} (\pa_{j_1} \ldots
   \pa_{j_5}[\eta_i\pa_{j_6}\phi])(\pa^{j_1} \ldots \pa^{j_5}[\eta_i\pa^{j_6}\phi])dx. \end{align}
    Now, \begin{align} \label{eq:may14080936} \pa^{j_1}
     \ldots \pa^{j_5} [\eta_i \pa^{j_6}\phi] = \eta_i (\pa^{j_1}
      \ldots \pa^{j_5} \pa^{j_6}\phi) + \sum (\na^{k_1} \eta_i)(\na^{k_2 + 1} \phi) \end{align}
       where the sum is over $k_1 + k_2 = 5$ such that $k_2 \leq 4$.
       To control the second term in (\ref{eq:may14080936}) we use the
        following procedure: Let $i_1 = i$. Suppose that we have found $i_1, \ldots, i_l$.
         The support of $\na^k \eta_{i_l}$ is contained in the image under $x$ of the
         set $\{y \in \R^2: 1 - \frac{d_0}{i_l} < |y| < 1 - \frac{d_0}{2i_l}\}$.
         Pick $i_{l+1}$ such that $1 - \frac{d_0}{2i_l} \leq 1 - \frac{d_0}{i_{l+1}}$.
Then $\eta_{i_{l+1}}$ takes the value $1$ on the set $\{y \in \R^2:
|y| \leq 1-\frac{d_0}{i_{l+1}}\}$ and $\{y \in \R^2: 1 -
\frac{d_0}{i_l} < |y| < 1-\frac{d_0}{2i_l}\} \subseteq \{y \in \R^2:
|y| \leq 1-\frac{d_0}{i_{l+1}}\}$.
\begin{lemma} \label{april29081455} For $k_2 \geq 1$ we have
\begin{align} (\na^{k_1} \eta_{i_1})(\na^{k_2}\phi) & =
\sum (\na^{k_1} \eta_{i_1})(\na^{l_2} \eta_{i_2})
\ldots (\na^{l_{n-1}} \eta_{i_{n-1}})(\na^{l_{n}}[\eta_{i_{n}}\na \phi]) \end{align}
 where the sum is over all $l_2 + \ldots + l_n = k_2 - 1$; where for instance
 if $l_2 = 0$ the term $\na^{l_2}\eta_{i_2}$ is taken to not be present in the sum;
 and where if $l_n = 0$ the term $\na^{l_{n}}[\eta_{i_n}\na \phi]$ is taken
 to be $\eta_{i_n}\na \phi$. \end{lemma}

\bp We prove this by induction on $k_2$. For $k_2 = 1$ we have
$(\na^{k_1} \eta_{i_1})(\na \phi) = \eta_{i_2}(\na^{k_1} \eta_{i_1})(\na \phi)$,
which is of the correct form. Suppose that $k_2 \geq 2$ and that we have the result
for smaller $k_2$. Then
\begin{align} (\na^{k_1} \eta_{i_1})(\na^{k_2}\phi) & = \eta_{i_2} \,(\na^{k_1}
\eta_{i_1})(\na^{k_2}\phi) \\ \label{eq:april29081402} &
= (\na^{k_1} \eta_{i_1})(\na^{k_2  - 1}[\eta_{i_2} \na \phi])
- (\na^{k_1} \eta_{i_1})\sum_{l_1 + l_2 = k_2 - 1,\, l_2 \leq k_2 - 2}
(\na^{l_1}\eta_{i_2})(\na^{l_2 + 1} \phi) \\ & = (\na^{k_1}
\eta_{i_1})(\na^{k_2 - 1}[\eta_{i_2} \na \phi]) \\
& - (\na^{k_1} \eta_{i_1})\sum_{l_1 + l_2 = k_2 - 1,\, l_2 \leq k_2 - 2}
\sum (\na^{l_1} \eta_{i_2})(\na^{m_2} \eta_{i_3}) \ldots (\na^{m_{n-1}}
\eta_{i_{n-1}})(\na^{m_{n}}[\eta_{i_{n}}\na \phi]) \end{align}
which again is of the correct form. \ep \bigskip

Integrating the first term in (\ref{eq:may14080936}) by parts twice we have
 \begin{align} & - \int_{\Om_t} (\pa_{j_1}
 \ldots \pa_{j_5} \pa^{j_6} [\eta_i \pa_{j_6}\phi])\eta_i(\pa^{j_1}
  \ldots \pa^{j_5} \phi)dx - \int_{\Om_t} (\pa_{j_1} \ldots \pa_{j_5}
   [\eta_i \pa_{j_6}f])( \pa^{j_6} \eta_i)(\pa^{j_1} \ldots \pa^{j_5} \phi)dx \\
   & \label{eq:may14080939} = \int_{\Om_t} (\pa_{j_1} \ldots \pa_{j_4} \pa^{j_6}
    [\eta_i \pa_{j_6}\phi])\eta_i(\pa^{j_1} \ldots \pa^{j_5} \pa_{j_5} \phi)dx
    + \int_{\Om_t} (\pa_{j_1} \ldots \pa_{j_4} \pa^{j_6}
    [\eta_i \pa_{j_6}\phi])(\pa_{j_5} \eta_i) (\pa^{j_1} \ldots \pa^{j_5} \phi)dx \\
    & - \int_{\Om_t} (\pa_{j_1} \ldots \pa_{j_5}
    [\eta_i \pa_{j_6}\phi])( \pa^{j_6} \eta_i)(\pa^{j_1} \ldots \pa^{j_5} \phi)dx \end{align}
  where we can control the second and third term in (\ref{eq:may14080939})
   using lemma \ref{april29081455}. Also,
    $\pa^{j_6} [\eta_i \pa_{j_6}\phi] = (\pa^{j_6}\eta_i)( \pa_{j_6}\phi) + \eta_i$
    and therefore the first term in (\ref{eq:may14080939}) is equal to
    \begin{align} \label{eq:may13080816} & \int_{\Om_t} \pa_{j_1} \ldots \pa_{j_4}
    [(\pa^{j_6}\eta_i)( \pa_{j_6}\phi)]\eta_i(\pa^{j_1} \ldots \pa^{j_5} \pa_{j_5} \phi)dx
    + \int_{\Om_t} \pa_{j_1} \ldots \pa_{j_4} [\eta_i]\eta_i(\pa^{j_1} \ldots
     \pa^{j_5} \pa_{j_5} \phi)dx \end{align} The above terms in (\ref{eq:may13080816})
     can be controlled using lemma \ref{april29081455} and the inductive hypothesis.
     This concludes the proof of proposition \ref{april29081232}.

\subsubsection{Boundary regularity.}

Let $\zeta$ denote $\zeta_1$. We have
$\na^2 [\zeta \phi] = (\na^2 \zeta)\,\phi + (\na \zeta)(\na \phi) + \zeta\,(\na^2 \phi)$
and $\na \zeta$ is supported in $\{y \in \R^2: 1-d_0 < |y| < 1 - \frac{d_0}{2}\}$.
We can find $i$ such that $1 - \frac{d_0}{2} \leq 1 - \frac{d_0}{i}$ and therefore
$\{y \in \R^2: 1-d_0 < |y| < 1 - \frac{d_0}{2}\}
\subseteq \{y \in \R^2: |y| < 1 - \frac{d_0}{i}\}$
where $\eta_i$ takes the value $1$.
Thus $(\na \zeta)(\na \phi) = \eta_i(\na \zeta)(\na \phi)
= (\na \zeta)(\na [\eta_i\phi]) - (\na \zeta)(\na \eta_i)\,\phi$,
which we control by theorem \ref{april29081232}.
Thus we need only be concerned with $\zeta \na^s \phi$.
In this section we prove the following result:

\begin{theorem} \label{april29081548} For $1 \leq s \leq 6$
\begin{align} \label{eq:april29081549}
\|\zeta \na^s \phi\|_{L^2(\Om_t)}\leq P\big[\|x\|_{5.5}\big] \end{align}
where $P$ is a polynomial.\end{theorem}

Since integration by parts on $\Om_t$ will yield a boundary term
which is difficult to deal with because $\pOm_t$ is the complement
of the singular support of $\phi$, we begin by expanding the region
of integration: Define $\tilde x = E(x)$ where $E$ is the extension
operator on $\Om$ $-$ see, for instance, \cite{MR0350177} $-$ and
define $\tilde V = \pa_t \tilde x$. Then both $\tilde V$ and $\tilde
x$ are defined in all of $\R^2$ and such that $\|\tilde
x\|_{H^s(\R^2)} \leq c\|x\|_{H^s(\Om)}$ and similarly for $\tilde
V$. Define $\tilde B^i_a = \frac{\pa \tilde x^i}{\pa y^a}$. Then
since $\det(B) = 1$ on $\Om$, we can choose $d_0$ (possibly smaller
than the one used before) so that $\tilde x$ is a change of
variables on $\tilde \Om = \{y \in \R^2: |y| < 1+ d_0\}$ and such
that $d_0$ is small enough that the normal $N$ to $\pOm_t$ can be
extended into the region between $\pOm_t$ and the boundary of
$\tilde \Om_t = \tilde x(t, \tilde \Om)$. This means that for every
$i$, $N$ can be extended into the support of $\zeta$ on both sides
of $\pOm_t$. Let $\tilde A$ be the inverse of $\tilde B$. We now
define $\tilde \phi$ as follows: \beq \label{eq:mar25081322} \tilde
\phi(t,x) = -  \chi_{\Om_t} * \Phi(x) \mbox{ for $x$ in } \tilde
\Om_t \eeq where again $\Phi$ is the fundamental solution for the
Laplacean. This means that on $\Om_t$, we have $\tilde \phi = \phi$
and therefore that $\tilde \phi$ and $\phi$ have the same regularity
on $\Om_t$. It also means that $\tilde \phi$ is smooth on $\pa
\tilde \Om_t$. Finally, let the norms on the extended domains
$\tilde \Om$ and $\tilde \Om_t$ be defined analogously to the norms
on $\Om$ and $\Om_t$.

Having extended the domain we now approximate $\tilde \phi$:
Let $\chi_m$ denote a smooth {\it radial} function compactly supported
in $\{y \in \R^2: |y| < 1 + \frac{1}{m}\}$, which takes the value $1$
on the set $\{y \in \R^2: |y| \leq 1-\frac{1}{m}\}$. This means that
$\pt \chi_m = 0$. By an abuse of notation we will also let $\chi_m$
denote the analogous function in the Eulerian frame. Define
$\tilde \phi_m(t,x) = -  \chi_m * \Phi(x) \mbox{ for $x$ in } \tilde \Om_t$. W
e now show that the approximations converge to $\tilde \phi$.

\begin{lemma} \label{mar31080845} $\|\na \tilde \phi_m -
\na \tilde \phi\|_{L^2(\tilde \Om_t)} \leq c \|\chi_m - \chi_{\Om_t}\|_{L^2(\tilde \Om_t)}$.
\end{lemma}

\bp From (\ref{eq:mar25081322}) it is clear that $\tilde \phi$ is in $H^1(\tilde \Om_t)$
 so integration by parts is justified. Similarly, for $\tilde \phi_m$.
 Now, \begin{align} \|\na \tilde \phi_m - \na \tilde \phi\|_{L^2(\tilde \Om_t)}^2
 & = \int_{\tilde \Om_t} (\pa_j \tilde \phi_m - \pa_j \tilde \phi)(\pa^j \tilde
 \phi_m - \pa^j \tilde \phi)dx \\ \label{eq:mar31080825} &
 = \int_{\pa \tilde \Om_t} N^j(\pa_j \tilde \phi_m -
 \pa_j \tilde \phi)(\tilde \phi_m - \tilde \phi)dS(x)
  - \int_{\tilde \Om_t} (\chi_m - \chi_{\Om_t})(\tilde\phi_m - \tilde\phi)dx. \end{align}
   To control the first term in (\ref{eq:mar31080825}) we note that there is $\delta>0$
such that $\mbox{dist}(\pa \tilde \Om_t, \pOm_t) > \delta$. This
means that for all $x$ on $\pa \tilde \Om_t$ and for all $z$ in
$\Om_t$ we have $|x-z| > \delta$. Hence for $x$ on $\pa \tilde
\Om_t$, \[|\na \tilde \phi_m(x) - \na \tilde \phi(x)| \leq
\int_{\tilde \Om_t} |\chi_m(z) - \chi_{\Om_t}(z)|\Phi'(|x-z|)dz \leq
\|\chi_m - \chi_{\Om_t}\|_{L^2(\tilde \Om_t)}. \] Also for $x$ in
$\pa \tilde \Om_t$ and for $x$ in $\tilde \Om_t$, $|\tilde\phi_m(x)
- \tilde\phi(x)| \leq \|\chi_m - \chi_{\Om_t}\|_{L^2(\tilde
\Om_t)}\|\Phi(|x-\cdot|)\|_{L^2(\tilde \Om_t)} \leq c\|\chi_m -
\chi_{\Om_t}\|_{L^2(\tilde \Om_t)}$, so we can control the first and
second term in (\ref{eq:mar31080825}) appropriately. \ep \bigskip

Lemma \ref{mar31080845} also shows that
$\|\zeta \na \phi_m - \zeta \na \phi\|_{L^2(\tilde \Om_t)} \leq P\big[\|x\|_2\big]
\|\chi_m - \chi_{\Om_t}\|_{L^2(\tilde \Om_t)}$. Let
$\tilde \phi_{m,n} = \tilde \phi_m - \tilde \phi_n$ and let $\chi_{m,n} = \chi_m - \chi_n$.
We will now prove the following proposition which shows that $(\zeta \na \pt^j \na \phi_m)$
is a Cauchy sequence in $L^2(\tilde \Om_t)$.

\begin{proposition} \label{april30082003} For $0 \leq j \leq 4$, we have
\begin{align} \|\zeta \na \pt^j \na \phi_{m,n}\|_{L^2(\Om_t)}
\leq P\big[\|x\|_{5.5}\big]\|\chi_{m,n}\|_{L^4(\tilde \Om_t)}, \end{align}
\end{proposition}

We begin by proving a lemma which says that we need only be concerned with tangential derivatives:

\begin{lemma} \label{philemma2a} Let $f$ satisfy $\Delta f = g$ in $\tilde \Om_t$
where $\pt g = 0$. Then for $0 \leq j \leq 4$,
\beq \label{eq:l22a} \|\zeta \na \pt^j \na f\|_{L^2(\tilde \Om_t)}
\leq P\big[\|x\|_{5.5}\big]\sum_{k = 0}^{j + 1}
\|\zeta \pt^k \na f\|_{L^2(\tilde \Om_t)} + \|g\|_{L^2(\tilde \Om_t)} \eeq
and for $0 \leq j \leq 2$ we have
\beq \label{eq:linfty2a} \|\zeta \na \pt^j \na f\|_{L^4(\tilde \Om_t)}
\leq P\big[\|x\|_{5.5}\big]\sum_{k = 0}^{j + 2}
\|\zeta \pt^k \na f\|_{L^2(\tilde \Om_t)} + \|g\|_{L^4(\tilde \Om_t)}.\eeq \end{lemma}

\bp We prove this result by induction. For $j = 0$ we have
$\|\zeta
\na^2 f\|_{L^2(\tilde \Om_t)} \leq \|\zeta \Delta f\|_{L^2(\tilde
\Om_t)} + \|\zeta \pt \na f\|_{L^2(\tilde \Om_t)}$, by lemma
\ref{normalderivnoninv}, which is of the right form. Now suppose
that $1 \leq j \leq 2$ and that we have (\ref{eq:l22a}) for smaller
$j$. Then $\|\zeta \na \pt^j \na f\|_{L^2(\tilde \Om_t)} \leq
\|\zeta \div \pt^j \na f\|_{L^2(\tilde \Om_t)} + \|\zeta \curl \pt^j
\na f\|_{L^2(\tilde \Om_t)} + \|\zeta \pt^{j+1} \na f\|_{L^2(\tilde
\Om_t)}$. Now $\zeta \div \pt^j \na f = \zeta \pt^j \Delta f + \sum
(\pt^k \tilde A)(\zeta \na \pt^l \na f)$ where the sum is over $k+l
= j$ such that $l \leq j-1 \leq 1$. Since $\pt^j \Delta f = 0$, we
have $\|\zeta \div \pt^j \na f\|_{L^2(\tilde \Om_t)} \leq
P\big[\|x\|_5]\|\zeta \na \pt^l \na f\|_{L^2(\tilde \Om_t)}$ which
we control by induction. Similarly, we also control $\|\zeta \curl
\pt^j \na f\|_{L^2(\tilde \Om_t)}$. For $j = 0$ we have
\begin{align} \|\zeta \na \pt^j \na f\|_{L^4(\tilde \Om_t)} & \leq
\|g\|_{L^4(\tilde \Om_t)} + \|\zeta \pt \na f\|_{L^4(\tilde \Om_t)}
\\ & \leq \|g\|_{L^4(\tilde \Om_t)} + \|\zeta \pt \na
f\|_{H^1(\tilde \Om_t)} \\ & \leq
P\big[\|x\|_{5.5}\big]\left[\|g\|_{L^4(\tilde \Om_t)} + \|\zeta
\pt^2 \na f\|_{L^2(\tilde \Om_t)}\right] \end{align} using Sobolev's
inequality and theorem \ref{april29081232}. Thus we have
(\ref{eq:linfty2a}) for $j=0$. Now suppose that $j=3$. Then for $0
\leq k \leq 2$ we have $\|(\pt^k \tilde A)(\zeta \na \pt^l \na
f)\|_{L^2(\tilde \Om_t)} \leq \|x\|_5\|\zeta\na \pt^l \na
f\|_{L^2(\tilde \Om_t)}$ and for $k = 3$ we have $l=0$ and $\|(\pt^k
\tilde A)(\zeta\na \pt^l \na f)\|_{L^2(\tilde \Om_t)} \leq \|\pt^3
\tilde A\|_{L^4(\tilde \Om_t)}\|\zeta\na^2 f\|_{L^4(\tilde \Om_t)}$,
both of which we control appropriately. Now we can prove
(\ref{eq:linfty2a}) for $j=1$ and using that result we can prove
(\ref{eq:l22a}) for $j=4$. And using that result we prove
(\ref{eq:linfty2a}) for $j=2$. \ep \bigskip

Using lemma \ref{philemma2a}, we see that it is enough to control
$\|\zeta \pt^j \na \tilde \phi_{m,n}\|_{L^2(\tilde \Om_t)}$ appropriately for
$0 \leq j \leq 5$ which is the content of the following proposition:

\begin{proposition} \label{mar31081044} For $0 \leq j \leq 5$ we have
\begin{align} \label{eq:april251014} \|\zeta
\pt^j \na \tilde \phi_{m,n}\|_{L^2(\tilde \Om_t)} \leq P\big[\|x\|_{5.5}\big]
\|\chi_{m,n}\|_{L^4(\tilde \Om_t)}. \end{align} \end{proposition}

\bp We prove that (\ref{eq:april251014}) holds by induction on the order.
To start the induction we have the analogue of lemma \ref{mar31080845}.
Now we suppose that we have $j = 5$ and suppose that we have already have
appropriate control of the lower order cases of (\ref{eq:april251014}).
We have \begin{align} \|\zeta \pt^5 \na \tilde \phi_{m,n}\|^2_{L^2(\tilde \Om_t)}
& = \int_{\tilde \Om_t} (\zeta \pt^5  \pa_i \tilde \phi_{m,n})
(\zeta \pt^5  \pa^i \tilde \phi_{m,n})dx \\ \label{eq:mar26081000}
& = \int_{\tilde \Om_t} (\zeta \pt^5 \pa_i
\tilde \phi_{m,n})(\zeta \pa^i \pt^5 \tilde \phi_{m,n})dx - \int_{\tilde \Om_t}
(\zeta \pt^5  \pa_i \tilde \phi_{m,n})(\pa^i \pt^5  x^l)(\zeta \pa_l \tilde \phi_{m,n})dx \\
& + \sum \int_{\tilde \Om_t} (\zeta \pt^5 \na \tilde \phi_{m,n})
(\na \pt^{l_1}x) \ldots (\na \pt^{l_{s-1}}x)(\zeta \pt^{l_s} \na \tilde \phi_{m,n})dx \end{align}
where $l_1 + \ldots + l_{s} = 5$ and $l_1, \ldots, l_{s}\leq 4$. For $l_1, \ldots, l_{s-1} \leq 2$
we control the third term in (\ref{eq:mar26081000}) by
$P\big[\|x\|_5\big]\|\zeta \pt^5 \na \tilde \phi_{m,n}\|_{L^2(\tilde \Om_t)}
\|\zeta \pt^{l_s} \na \tilde \phi_{m,n}\|_{L^2(\tilde \Om_t)}$,
which we control by induction. Suppose that $3 \leq l_1 \leq 4$.
Then $\|\na \pt^{l_1} x\|_{L^4(\tilde \Om_t)} \leq \|x\|_{5.5}$;
and $l_2, \ldots, l_{s - 1} \leq 2$, so we control the other terms containing $x$.
We also have $0 \leq l_s \leq 2$ and therefore $\|\zeta \pt^{l_s} \na
\tilde \phi_{m,n}\|_{L^4(\tilde \Om_t)} \leq \|\zeta \na \pt^{l_s} \na
\tilde \phi_{m,n}\|_{L^2(\tilde \Om_t)}$, which we control appropriately by lemma
\ref{philemma2a}. Thus we control the third term in (\ref{eq:mar26081000}).
Integrating the first two terms in (\ref{eq:mar26081000}) by parts gives
\begin{align} \label{eq:mar27080822} & \int_{\tilde \Om_t} (\zeta \pt^5
\pa_i \tilde \phi_{m,n})(\pt^5 x^l)(\zeta \pa^i \pa_l \tilde \phi_{m,n})dx
- \sum \int_{\tilde \Om_t} (\zeta \pa^i \pt^5 \pa_i \tilde \phi_{m,n})
(\pt^{j_1} x)(\zeta \pt^{j_2} \na \tilde \phi_{m,n})dx \\ & + \sum \int_{\tilde \Om_t}
(N^i \pt^5  \pa_i \tilde \phi_{m,n})(\pt^{j_1}x)(\pt^{j_2} \na \tilde \phi_{m,n})dS(x) \end{align}
where the sums are over all $j_1 + j_2 = 5$ such that $j_1, j_2 \leq 4$.
Here we are ignoring the terms which arise from the derivative falling on
$\zeta$ because in this case we can use theorem \ref{april29081232}.
We control the first term in (\ref{eq:mar27080822}). Also,
\beq \label{eq:mar26081059} \zeta \pa^i \pt^5  \pa_i \tilde \phi_{m,n}
= (\na \pt^5  x)(\zeta \na^2 \tilde \phi_{m,n}) + \sum (\na \pt^{l_1}x)
\ldots (\na \pt^{l_{s-1}}x)(\zeta \na \pt^{l_s} \na \tilde \phi_{m,n}) \eeq
where $l_1 + \ldots + l_{s} = 5$ and $l_1, \ldots, l_{s} \leq 4$, since
$\pt^5 \Delta \tilde \phi_{m,n} = 0$. For $l_1, \ldots, l_{s-1} \leq 2$
we control the second term in (\ref{eq:mar26081059}) appropriately by
lemma \ref{philemma2a}. If $3 \leq l_1 \leq 4$ then we control
$\na \pt^{l_1} x$ as before and $0 \leq l_2, \ldots, l_{s-1} \leq 2$ so we can control
the other terms containing $x$ in $L^\infty(\tilde \Om_t)$. We also have $0 \leq l_s \leq 2$
so we can control $\|\zeta \na \pt^{l_s} \na \tilde \phi_{m,n}\|_{L^4(\tilde \Om_t)}$ by
lemma \ref{philemma2a}. We therefore control the second term in (\ref{eq:mar26081059}).
Let us now consider the first term in (\ref{eq:mar26081059}): Commute one $\pt$ to
the outside to obtain, in addition to a lower order term,
\begin{align} \int_{\tilde \Om_t} (\pt \na \pt^4  x)(\zeta \na^2
\tilde \phi_{m,n})(\pt^{j_1} x)(\zeta \pt^{j_2} \na \tilde \phi_{m,n})dx
\label{eq:mar26081128} & = - \int_{\tilde \Om_t} (\na \pt^4 x)(\zeta \pt \na^2
\tilde \phi_{m,n})(\pt^{j_1} x)(\zeta \pt^{j_2} \na \tilde \phi_{m,n})dx \\ &
- \int_{\tilde \Om_t} (\na \pt^4 x)(\zeta \na^2 \tilde \phi_{m,n})(\pt^{j_1+1} x)(
\zeta \pt^{j_2} \na \tilde \phi_{m,n})dx \\ & - \int_{\tilde \Om_t} (\na \pt^4 x)(
\zeta \na^2 \tilde \phi_{m,n})(\pt^{j_1}  x)(\zeta \pt^{j_2+1} \na \tilde \phi_{m,n})dx \end{align}
where no boundary terms arise because the components of $\pt$ are orthogonal to the normal
on $\pa \tilde \Om_t$. In all of the terms in (\ref{eq:mar26081128}) we control the first
two factors in each integrand using lemma \ref{philemma2a}. In the second term in
(\ref{eq:mar26081128}) we also have, for $j_1 \leq 2$,
$\|(\pt^{j_1+1}  x)(\zeta \pt^{j_2} \na \tilde \phi_{m,n})\|_{L^2(\tilde \Om_t)}
\leq \|x\|_5\|\zeta \pt^{j_2} \na \tilde \phi_{m,n}\|_{L^2(\tilde \Om_t)}$ which we control.
For $3 \leq j_1 \leq 4$ we have $1 \leq j_2 \leq 2$ and therefore we control the second
term in (\ref{eq:mar26081128}) in this case also. The third term in (\ref{eq:mar26081128})
follows similarly.

Now we control the boundary term in (\ref{eq:mar27080822}): \beq
\label{eq:mar26081213} \sum \int_{\pa \tilde \Om_t} (\pt^{k_1} x)
\ldots (\pt^{k_s} x) (\pa_i \pa_{k_1} \ldots \pa_{k_s}\tilde
\phi_{m,n})N^i(\pt^{l_1} x) \ldots (\pt^{l_s} x) (\pa_{l_1} \ldots
\pa_{l_s}\tilde \phi_{m,n})dS(x) \eeq where the sum is over $k_1 +
\ldots + k_s = 5$ and $l_1 + \ldots + l_s = 5$ such that $l_1,
\ldots, l_s \leq 4$. As was mentioned above, there is $\delta>0$
such that for all $x$ on $\pa \tilde \Om_t$ and $z$ in $\Om_t$ we
have $|x-z|>\delta$. Therefore $|\na^s \tilde \phi_{m,n}(x)| \leq
\|\chi_{m,n}\|_{L^2(\tilde \Om_t)}$. The highest order term of the
above terms is \beq \label{eq:mar26081223} \int_{\pa \tilde \Om_t}
(\pt^5x)(\na \tilde \phi_{m,n})N(\pt^4 x) (\pt x)(\na ^2\tilde
\phi_{m,n})dS(x) \eeq which is controlled by
$P\big[\|x\|_{5.5}\big]$ using the trace theorem. This
concludes the proof. \ep \bigskip

By lemma \ref{philemma2a}, therefore, we have proposition \ref{april30082003},
which means that $(\zeta \na \pt^j \na \tilde \phi_m)_{m = 1}^\infty$ is a
Cauchy sequence in $L^2(\tilde \Om_t)$, for $0 \leq j \leq 4$. This means that
$\zeta \na \pt^j \na \tilde \phi_m \to \zeta \na \pt^j \na \tilde \phi$ in
$L^2(\tilde \Om_t)$. From lemma \ref{philemma2a} and proposition \ref{mar31081044}
we have $\|\zeta \na \pt^j \na \tilde \phi_m\|_{L^2(\tilde \Om_t)}
\leq P\big[\|x\|_{5.5}\big]$, for $0 \leq j \leq 4$ and therefore have
$\|\zeta \na \pt^j \na \tilde \phi\|_{L^2(\tilde \Om_t)}\leq P\big[\|x\|_{5.5}\big]$
and hence \begin{align} \label{eq:april30082056} \|\na \pt^j \na \phi\|_{L^2(\Om_t)}
\leq P\big[\|x\|_{5.5}\big], \end{align} for $0 \leq j \leq 4$. From now on we
no longer consider the extended domain; all norms are now the usual, non-extended, norms.

\begin{lemma} \label{philemma2b} Let $f$ satisfy $\Delta f = g$ in $\Om_t$ where
$\pt g = \pa_r g = 0$ on $\Om_t$. Writing $\pa$ to denote both $\pt$ and $\pa_r$
we have \beq \label{eq:l22ab} \|\zeta \na \pa^j \na f\|_{L^2(\Om_t)}
\leq P\big[\|x\|_{5}\big]\sum_{k = 0}^{j} \|\zeta \pt \pa^k \na f\|_{L^2(\Om_t)}
+ \|g\|_{L^2(\Om_t)} \eeq for $0 \leq j \leq 4$; and for $0 \leq j \leq 2$ we have
\beq \label{eq:linfty2ab} \|\zeta \na \pa^j \na f\|_{L^\infty(\Om_t)}
\leq P\big[\|x\|_{5}\big]\sum_{k = 0}^{j + 1} \|\zeta \pt^2 \pa^k \na f\|_{L^2(\Om_t)}
+ \|g\|_{L^\infty(\Om_t)}.\eeq \end{lemma}

\bp We prove this result by induction. For $j = 0$ we have
$\|\zeta \na^2 f\|_{L^2(\Om_t)} \leq \|\zeta \Delta f\|_{L^2(\Om_t)}
+ \|\zeta \pt \na f\|_{L^2(\Om_t)}$, by lemma \ref{normalderivnoninv},
which is of the right form. Now suppose that $1 \leq j \leq 2$ and that we
have (\ref{eq:l22a}) for smaller $j$. Then
$|\zeta \na \pa^j \na f| \leq |\zeta \div \pa^j \na f| + |\zeta \curl \pa^j \na f|
+ |\zeta \pt \pa^{j} \na f|$. Now $\zeta \div \pa^j \na f = \zeta \pa^j \Delta f
+ \sum (\pa^k A)(\zeta \na \pa^l \na f)$ where the sum is over $k+l = j$ such
that $l \leq j-1 \leq 1$. Since $\pa^j \Delta f = 0$, we have
$\|\zeta \div \pa^j \na f\|_{L^2(\Om_t)} \leq P\big[\|x\|_5]\|\zeta \na \pa^l \na f\|_{L^2(\Om_t)}$
which we control by induction. Similarly, we also control $\|\zeta \curl \pa^j \na f\|_{L^2(\Om_t)}$.
For $j = 0$ we have
\begin{align} \|\zeta \na \pa^j \na f\|_{L^\infty(\Om_t)} & \leq \|g\|_{L^\infty(\Om_t)}
+ \|\zeta \pt \na f\|_{L^\infty(\Om_t)} \\   & \leq \|g\|_{L^\infty(\Om_t)}
+ \|\zeta \na \pa \pt  \na f\|_{L^2(\Om_t)} \\
& \leq P\big[\|x\|_{5}\big]\left[\|g\|_{L^\infty(\Om_t)}
+ \|\zeta \pa \pt^2 \na f\|_{L^2(\Om_t)}\right] \end{align} using
Sobolev's inequality and lemma \ref{philemma2a}. Thus we have (\ref{eq:linfty2a})
for $j=0$. Now suppose that $j=3$. Then for $0 \leq k \leq 2$ we have
$\|(\pa^k A)(\zeta \na \pa^l \na f)\|_{L^2(\Om_t)} \leq \|x\|_5\|\zeta\na \pa^l \na f\|_{L^2(\Om_t)}$
and for $k = 3$ we have $l=0$ and
$\|(\pa^k A)(\zeta\na \pa^l \na f)\|_{L^2(\Om_t)} \leq \|\pa^3 A\|_{L^2(\Om_t)}
\|\zeta\na^2 f\|_{L^\infty(\Om_t)}$, both of which we control appropriately.
Now we can prove (\ref{eq:linfty2a}) for $j=1$ and using that result we can prove
(\ref{eq:l22a}) for $j=4$. And using that result we prove (\ref{eq:linfty2a})
for $j=2$. \ep \bigskip

Using lemma \ref{philemma2b} and an induction argument we control
$\|\zeta \na^s \phi\|_{L^2(\Om_t)}$ for $0 \leq s \leq 5$ and hence we obtain
theorem \ref{ChopinEtudeAFlat}. \ep

\subsection{Estimates for $p$.}

Taking the divergence of (\ref{eq:DebussyTroisPoemesMov1}) we see
that $p$ satisfies $\Delta p = -(\pa_i v^j)(\pa_j v^i) + 1$ on
$\Om_t$ and $p = 0$ on $\pOm_t$. Thus we have the following:

\begin{proposition} \label{mar25081252} For all $i \geq 1$ we have
\begin{align} \label{eq:may03081030} \|\eta_i \na p\|_{H^5(\Om_t)}
\leq P\big[\|x\|_{5}, \,\|V\|_{5}\big] \end{align} and
\begin{align} \label{eq:may03081031} \|\na p\|_{H^4(\Om_t)} \leq P\big[\|x\|_{5}, \,\|V\|_4\big] \end{align}
for $0 \leq j\leq 4$. Moreover, we have
\begin{align}\label{eq:may03081032}\|\na p\|_{H^{4.5}(\Om_t)}
\leq P\big[\|x\|_{5.5}, \,\|V\|_{4.5}\big] \end{align} and
\begin{align}\label{eq:may03081033}\|\na \dot p\|_{H^3(\Om_t)}
\leq P\big[\|x\|_{5.5}, \,\|V\|_{5}\big] \end{align}
where $\dot p = \pa_t p$. \end{proposition}

\bp The estimate (\ref{eq:may03081030}) follows similarly to theorem \ref{april29081232}.
The estimates (\ref{eq:may03081031}) and (\ref{eq:may03081032}) follow similarly to
theorem \ref{ChopinEtudeAFlat}. The estimate (\ref{eq:may03081033})
follows similarly to theorem \ref{ChopinEtudeAFlat} using theorem
\ref{ChopinEtudeAFlat} and (\ref{eq:may03081031}) above.
For a detailed explanation of these proofs, see \cite{mythesis}. \ep

\section{Energy estimates.}\label{mar26080821}

Finally, we are ready to prove the energy estimate in theorem
\ref{may16080841}. We control $\|V\|_5$ using lemma
\ref{april04080805} and lemma \ref{CSineq}, together with $E_1$ and
$E_2$ below: Let \beq \label{eq:mar25081109} E_1(t) =  \left\|\zeta
\pt^5 v\right\|_{L^2(\Om_t)}^2 + \left\|\eta
v\right\|_{H^5(\Om_t)}^2 \eeq and $E_2(t) = \|\curl
(v)\|_{H^{4.5}(\Om_t)}$, where $\zeta = \zeta_1$ is the cut-off
function supported near the boundary of $\Om_t$ and $\eta = \eta_1$
is the cut-off function supported in the interior of $\Om_t$, as
defined in section \ref{april02081714}. Note that we will ignore
terms which arise from the derivative falling on the cut-off
function because these terms will be of lower order. To build
regularity for $\|x\|_{5.5}$ we use lemma \ref{CSineq} together
with $E_3$ and $E_4$ below: Let \beq \label{eq:mar31080739} E_3(t) =
\int_{\pOm_t} (- \na p \cdot N) \left[(\pt^5 x) \cdot N\right]^2
dS(x), \eeq where $N$ is the external unit normal to $\pOm_t$. Note
that the term in (\ref{eq:mar31080739}) should read $\int_{\pOm_t}
(- \na p \cdot N) \left[(\pt^5 x) \circ x^{-1} \cdot N\right]^2 \\
dS(x)$, but in the interest of creating tidy computations we will
hereinafter omit some terms $-$ such as the `$\circ x^{-1}$' above
$-$ which are not crucial to understanding the argument. And let \[
E_4(t) = \left\|\div \left[\pa x\right]
\right\|_{H^{3.5}(\Om_t)}^2 + \left\|\curl \left[\pa x\right]
\right\|_{H^{3.5}(\Om_t)}^2 \] where $\div[\pa x] = \pa_i [(\pa
x^i)\circ x^{-1}]$ and $\pa$ is an arbitrary derivative in the
Lagrangian frame. We now explain how $E_3$ and $E_4$ will be used to
bound $\|x\|_{5.5}$. Let $\sum x_k e^{ik\theta}$ be the
tangential Fourier expansion of $x$. Then $(\ptt^5 x) \cdot N = \sum
\langle k \rangle^5 x_k \cdot N e^{ik\theta}$ and $(\pt^j x) \cdot N
= \sum (ik)^j x_k \cdot N e^{ik\theta}$. Therefore \begin{align}
\|(\ptt^5 x) \cdot N\|_{L^2(\pOm_t)}^2 \leq \sum \langle k
\rangle^{10} |x_k \cdot N|^2 \leq \sum |x_k \cdot N|^2 + \sum k^2
|x_k \cdot N|^2 \ldots + \sum k^{10} |x_k \cdot N|^2  \end{align}
and \begin{align} \sum k^{2j} |x_k \cdot N|^2 \leq \|(\pt^j x) \cdot
N\|_{L^2(\pOm_t)}^2. \end{align} Using the trace theorem and the
fact that $x(t,y) = y + \int_{[0,t]} V(s,y)ds$ we control the terms
with $0 \leq j \leq 4$. For the highest order term we have \beq
\label{eq:feb2008115} \left\|\left(\pt^5 x \right) \cdot N
\right\|_{L^2(\pOm_t)}^2 =  \int_{\pOm_t} \left[ \left(\pt^5 x
\right) \cdot N\right]^2 dS(x) \leq \frac{1}{c_0}\int_{\pOm_t} (-
\na p \cdot N)\left[\left(\pt^5 x \right) \cdot N\right]^2 dS(x)
\leq \frac{E_3}{c_0} . \eeq Using $E_4$ we control $\left\|\div
\left[\pa x\right] \right\|_{H^{3.5}(\Om_t)}^2 \mbox{ and }
\left\|\curl \left[\pa x\right] \right\|_{H^{3.5}(\Om_t)}^2$,
and therefore, from lemma \ref{CSineq}, we have, by
(\ref{eq:feb2008115}), \b \label{eq:feb18080910} \left\|\pa \pt
x\right\|^2_{H^{4.5}(\Om_t)} \leq P\big[\|x\|_5\big]\left[E_4 +
\frac{E_3}{c_0}\right]. \e This means that we similarly control
$\left\|\pt \pa x\right\|^2_{H^{3.5}(\Om_t)}$ and therefore
$\left\|\pa x\right\|^2_{H^{4}(\pOm_t)}$. Hence we control
$\left\|\pa x\right\|^2_{{4.5}}$ and therefore
$\|x\|_{5.5}$.

\subsection{Almost $E_1$.} \label{E1}

The time derivative of $E_1$ is equal to
\begin{align} \label{eq:feb16080806} - 2 \int_{\Om_t} (\zeta \pt^5 v^i)
(\zeta \pt^5 \pa_i p) dx - 2 \int_{\Om_t} (\zeta \pt^5 v^i)
(\zeta \pt^5 \pa_i \phi) dx - 2 \int_{\Om_t} (\eta \pt^5 v^i) (\zeta \pt^5 \pa_i p) dx
- 2 \int_{\Om_t} (\eta \pt^5 v^i) (\zeta \pt^5 \pa_i \phi) dx \end{align} using
(\ref{eq:DebussyTroisPoemesMov1}). Using proposition \ref{april29081232} and
(\ref{eq:may03081030}) from proposition \ref{mar25081252} we control the third and
fourth term in (\ref{eq:feb16080806}). The second term in (\ref{eq:feb16080806})
can be controlled using theorem \ref{april29081548}. It now remains to control the
first term in (\ref{eq:feb16080806}). We will deal with this term in section \ref{therest}.

\subsection{$E_2$.} \label{E2}

We have $[\pa_t, \pa_i]x^j = -(\pa_i v^j)$ and therefore
\begin{align} \pa_t \curl(v) & = \pa_1 \pa_t v_2 - \pa_2 \pa_t v_1
- (\pa_1 v^j)(\pa_j v_2) + (\pa_2 v^k)(\pa_k v_1) \\ & = - (\pa_1
v_1)(\pa_1 v_2)
 (\pa_1 v_2)(\pa_2 v_2) + (\pa_2 v_1)(\pa_1 v_1)  +
(\pa_2 v_2)(\pa_2 v_1) \\ & = - (\pa_1 v_2)[(\pa_1 v_1) + (\pa_2
v_2)] + (\pa_2 v_1)[(\pa_1 v_1) + (\pa_2 v_2)] \\ & =
-\curl(v)\div(v) \\ & = 0. \end{align} Thus $\curl(v)(t) =
\curl(v)(0)$ and therefore $\|\curl(v)(t)\|_{H^{4.5}(\Om_t)} =
\|\curl(v)(0)\|_{H^{4.5}(\Om_t)}$.

\subsection{$E_3$.} \label{E3}

The time derivative of $E_3$ is equal to
\begin{align} \label{eq:may16080853} & \int_{\pOm_t} \pa_t |\na p|
\left[(\pt^5 x) \cdot N\right]^2 dS(x) + 2\int_{\pOm_t} |\na p|
\left[(\pt^5 x) \cdot (\pa_t N) \right]\left[(\pt^5 x) \cdot N\right] dS(x) \\
& + 2 \int_{\pOm_t} |\na p| \left[(\pa_t \pt^5 x) \cdot N\right]\left[(\pt^5 x)
\cdot N\right]dS(x). \end{align} In (\ref{eq:may16080853}) we control the first
and second term. Using the fact that $\frac{-\pa_i p}{|\na p|} = N_i$ on the
third term in (\ref{eq:may16080853}) we have
\begin{align} \label{eq:may16080859} - 2 \int_{\pOm_t} (\pt^5 v^i)N_i(\pt^5 x^j)(\pa_j p)dS(x)
& = - 2 \int_{\Om_t} (\pa_i \pt^5 v^i)(\pt^5 x^j)(\pa_j p)dx
- 2 \int_{\Om_t} (\pt^5 v^i)(\pa_i \pt^5 x^j)(\pa_j p)dx \\
& - 2 \int_{\Om_t} (\pt^5 v^i)(\pt^5 x^j)(\pa_i \pa_j p)dx \end{align}
using the divergence theorem. We control the third term in (\ref{eq:may16080859}).
In the first term in (\ref{eq:may16080859}) we commute a $\pt$ outside the $\pa_i$
this generates a lower order term and also
\begin{align} \label{eq:april04081043} - 2 \int_{\Om_t} (\pt \pa_i \pt^4 v^i)
(\pt^5 x^j)(\pa_j p)dx \leq \|\ptt^\frac12 \pa_i \pt^4 v^i \|\|\ptt^\frac12[(\pt^5 x)(\na p)]\| \end{align}
using lemma \ref{HalfIBP}. By lemma \ref{algebra} and the fact that $\div[v] = 0$,
we can control the above. The second term in (\ref{eq:may16080859}) remains.
We deal with this term in the next section.

\subsubsection{The rest.} \label{therest}

Combining the first term from (\ref{eq:feb16080806}) and the second
term from (\ref{eq:may16080859}) gives \begin{align}
\label{eq:may03081056} & = - 2 \int_{\Om_t} (\zeta \pt^5v^i)(\zeta
\pa_i \pt^5 p)dx + 2 \sum \int_{\Om_t} \zeta^2 (\pt^5 v)(\na
\pt^{k_1} x) \ldots (\na \pt^{k_{s-1}} x)(\pt^{k_{s}}\na p)
\end{align} where $k_1 + \ldots + k_s = 5$ and $k_1, \ldots, k_s
\leq 4$. In the first term in (\ref{eq:may03081056}) we integrate by
parts to obtain \begin{align} \int_{\Om_t} (\zeta \pa_i
\pt^5v^i)(\zeta \pt^5 p)dx. \end{align} Again, we integrate half of
one of the $\pt$ from $\div \pt^5 v$ to the other side. The result
can be controlled by (\ref{eq:may03081032}) from proposition
\ref{mar25081252} and an argument from above. We can control the sum
in (\ref{eq:may03081056}) using (\ref{eq:may03081031}) from
proposition \ref{mar25081252}.

\subsection{$E_4$.} \label{E4}

First we deal with the divergence term. We have \beq
\label{eq:feb20081337} \pa_t \div [\pa x] = (\na v)(\pa^2 x) + \div
[\pa_t \pa x] = (\na v)(\pa^2 x) + \na \div v. \eeq Therefore we
have an equation of the form $\pa_t f = g$. Since
$H^{3.5}(\Om_t)$ is an algebra, we control the first term in
(\ref{eq:feb20081337}) by $\|V\|_{4.5}\|x\|_{5.5}$. We now
consider two time derivatives on $\curl [\pa x]$: \b \pa_t^2 \curl
[\pa x] & = & \pa_t \left[(\na v)(\pa^2 x) + \curl \pa_t [\pa
x]\right] \\ & = & \pa_t \left[(\na v)(\pa^2 x)\right] + (\na
v)(\na^2 v) + \curl \pa_t^2 [\pa x] \\ \label{eq:feb22081011}  & = &
\pa_t \left[(\na v)(\pa^2 x)\right] - (\na \pa_t v)(\pa^2 x) + (\pa
x)(\pa \pa_t^2 x) \e since $\curl \pa_t^2 x = 0$. Equation
(\ref{eq:feb22081011}) is of the form $\pa_t[(\pa_tf) - g] = h$.
Integrating with respect to time once yields $(\pa_t f)(t) =
(\pa_tf)(0) + g(t) - g(0) + \int_{[0,t]}h(u)du$. Another integration
with respect to time again gives \begin{align} f(t) = f(0) + t(\pa_t
f)(0)  + \int_{[0,t]} g(u)du - tg(0) +
\int_{[0,t]}\int_{[0,u_2]}h(u_1)du_1du_2.
\end{align} Here $f = \curl [\pa x]$ so we control $f(0)$ and
$(\pa_t f)(0)$. We have already seen that we can control the first
term in (\ref{eq:feb22081011}). The second term in
(\ref{eq:feb22081011}) can be controlled using the fact that
$H^{3.5}(\Om_t)$ is an algebra, (\ref{eq:may03081032}) from
proposition \ref{mar25081252} and theorem \ref{ChopinEtudeAFlat}.

\appendix

\section{Properties of $\ptt$.}

In this section we prove a result concerning how $\ptt^\frac12$ acts on a product and also
an integration by parts type result.

\begin{lemma} \label{algebra} Let $f$ and $g$ be functions on $\Om$. Then
$\|\ptt^\frac12[fg] - \ptt^\frac12[f]g\|^2 \leq c\|f\|^2\|\ptt^{\frac12 + a} g\|^2$
for $a>\frac12$. \end{lemma}

\bp Let $\sum f_k(\rho)e^{ik\theta}$ and $\sum g_l(\rho)e^{il\theta}$ be tangential
Fourier expansions of $f(\rho,\theta)$ and $g(\rho,\theta)$ respectively.
Then \[ fg = \sum_k \sum_l f_kg_l e^{i[k+l]\theta}
= \sum_m \left[\sum_{k+l = m} f_kg_l \right]e^{im\theta} \] and therefore
\beq \label{eq:mar27080843} \ptt^\frac12[fg] = \sum_m \langle m \rangle^\frac12
 \left[\sum_{k+l = m} f_kg_l \right]e^{im\theta} \eeq where
 $\langle m \rangle = [1 + |m|^2]^\frac12$. Also
 $\ptt^\frac12[f] = \sum_k \langle k \rangle^\frac12 f_ke^{ik\theta}$ and therefore
 \beq \label{eq:mar27080844} \ptt^\frac12[f]g = \sum_m \left[\sum_{k+l=m}
 \langle k \rangle^\frac12 f_kg_l\right]e^{im\theta}. \eeq The difference between
 (\ref{eq:mar27080843}) and (\ref{eq:mar27080844}) is
 \beq \label{eq:mar27080845} \sum_m \sum_{k+l=m} \left[\langle k+l \rangle^\frac12
 - \langle k \rangle^\frac12\right]f_kg_l e^{im\theta}. \eeq
 We can control this using the following lemma.

\begin{lemma} \label{xes} Let $k$ and $l$ be points in ${\bf Z}$.
Then \[\left|\langle k + l \rangle^\frac12 - \langle k \rangle^\frac12\right|
\leq c \langle  l \rangle^{\frac12}, \] where $c$ is a constant. \end{lemma}

\bp Suppose that $k$ and $l$ are such that $0 \leq |k| \leq |l|$.
Then \[\Big|\langle k+l\rangle^\frac12 - \langle k
\rangle^\frac12\Big| \leq c\langle l\rangle^\frac12. \] Now suppose
that $k$ and $l$ are such that $0 \leq |l| < |k|$. Then \bea \langle
k\rangle^\frac12 \left|\frac{(1 + (k+l)(k+l) )^{\frac{1}{4}}}{(1 +
k^2)^{\frac{1}{4}}} - 1\right| &=& \langle k\rangle^\frac12
\left|\left(\frac{1 + k^2 + 2k l + l^2}{1 +
k^2}\right)^{\frac{1}{4}} - 1\right| \\ &=& \langle k\rangle^\frac12
\left|\left(1 + \frac{2k l + l^2}{1 +  k  ^2}\right)^{\frac{1}{4}} -
1\right|.\eea Define $f(x) = (1 + x)^{\frac{1}{4}} - 1 $. Then there
is a constant $c$ which bounds $\frac{f(x)}{|x|}$, for all $x$ in
$(-4,4)$. Therefore, \begin{align} \langle k\rangle^\frac12
\left|\left(1 + \frac{2k  l + l^2}{1 + k^2}\right)^{\frac{1}{4}} -
1\right| \leq c\langle k\rangle^\frac12   \frac{2kl + l^2}{\langle
k\rangle^2} \leq c \langle k\rangle^\frac12 \frac{\langle k \rangle
\langle l \rangle + \langle l \rangle^2 }{\langle k\rangle^2} \leq
c\frac{ \langle l \rangle}{\langle k\rangle^{\frac12}}. \end{align}
Since $|l| < |k|$, $\frac{\langle k\rangle^{\frac12}}{ \langle l
\rangle^{\frac12}} > 1$, from where the result follows. \ep \bigskip

And from lemma \ref{xes} we see that (\ref{eq:mar27080845}) can be estimated in $L^2(\Om)$
by \begin{align} \sum_m \left[\sum_{k+l = m} \langle l \rangle^\frac12 \frac{\langle l
\rangle^a}{\langle l \rangle^a}|f_k||g_l| \right]^2 & \leq \sum_m \left[\sum_{k+l = m}
\langle l \rangle^{2(\frac12 + a)}|g_l|^2\right]\left[\sum_{k+l = m} \langle l
\rangle^{-2a} |f_k|^2\right]. \end{align} Since $\sum \langle l \rangle^{-2a}$
is convergent for $2a>1$ we must have $a>\frac12$. This proves lemma \ref{algebra}. \ep \bigskip

From this proof we also have the following result:

\begin{corollary} \label{algebra2} Let $f$ and $g$ be functions on $\pOm$.
Then \begin{align} \|\ptt^\frac12[fg] - \ptt^\frac12[f]g\|_{L^2(\pOm)}^2
\leq c\|f\|_{L^2(\pOm)}^2\|\ptt^{\frac12 + a} g\|_{L^2(\pOm)}^2 \end{align} for $a > \frac12$. \end{corollary}

\begin{lemma} \label{HalfIBP} Let $f$ and $g$ be functions on $\Om$.
and let $( \mbox{ }, \mbox{ })$ be the $L^2(\Om)$-inner product,
then \newline $|(f, \pt g)| \leq c\|\ptt^\frac12 f
\|\|\ptt^\frac12g\|$. \end{lemma}

\bp Let $\sum f_k(\rho)e^{ik\theta}$ and $\sum g_l(\rho)e^{il\theta}$ be
tangential Fourier expansions of $f(\rho,\theta)$ and $g(\rho,\theta)$ respectively.
Then \begin{align} |(f, \pt g)| & = \left|\int_0^1 \int_{0}^{2\pi} f(\rho,\theta)g(\rho,\theta)
\rho d\rho d \theta \right| \\ & \leq \left|\int_0^1 \sum_l \langle l \rangle^\frac12 f_l(\rho)\langle
 l \rangle^\frac12 g_l(\rho) \rho d\rho \right| \\ & \leq \|\ptt^\frac12f\|\|\ptt^\frac12g\|. \end{align}
  \ep

\section{Hodge-decomposition inequalities.}

In this section we prove the results whose proofs were omitted in the body of the text.
We begin with a lemma which says that in the support of $\zeta$ we can control all
derivatives by the curl the divergence and a tangential derivative.

\begin{lemma} \label{normalderivnoninvII} Let $\alpha$ be a vector-field on $\tilde \Om_t$.
Define $(\curl \alpha)_{jk}=\pa_j \alpha_k-\pa_k \alpha_j$ and $\div \alpha = \pa_j \alpha^j$.
Then we have the following pointwise estimate on $\Om_t$:
\beq |\zeta \na \alpha| \leq |\zeta \curl \alpha| + |\zeta \div \alpha| + |\zeta \,\pt \alpha|,\eeq
where $|\cdot|$ denotes the usual Euclidean distance. \end{lemma}

\bp Here we will suppress the index on $\zeta$, letting it be denoted simply by $\zeta$.
Define $(\de \alpha)_{jk} = \pa_j \alpha_k + \pa_k \alpha_j$. Thus
$2\na \alpha = \curl \alpha + \de \alpha$. Let
$\beta = \mbox{diag}(\pa_1 \alpha_1, \ldots, \pa_n \alpha_n)$ and define
$\gamma = \zeta \de \alpha - \zeta \beta$. Then $|\zeta \na \alpha| \leq |\zeta \curl \alpha|
+ |\zeta \div \alpha| + |\gamma|$. It remains to control $\gamma$. Also define
\beq \label{eq:jan09080916} {Q}^{jk} = \delta^{jk} - N^jN^k, \eeq the projection onto
tangential vector-fields. Hence
\bea |\gamma|^2 & = & \delta^{ij}\delta^{kl} \gamma_{ik}\gamma_{jl} \\ &
= & \left({Q}^{ij}+{N}^i {N}^j\right)\left({Q}^{kl}+{N}^k {N}^l\right)
\gamma_{ik}\gamma_{jl} \\ & = & {Q}^{ij}{Q}^{kl}\gamma_{ik}\gamma_{jl}
+ {Q}^{ij}{N}^k {N}^l \gamma_{ik}\gamma_{jl} + {N}^i {N}^j {Q}^{kl} \gamma_{ik}
\gamma_{jl} \\ & + &  {N}^i {N}^j {N}^k {N}^l\gamma_{ik}\gamma_{jl}.\eea Since
$\gamma$ is symmetric, ${N}^i {N}^j {Q}^{kl}\gamma_{ik}\gamma_{jl} = {Q}^{ij} {N}^k {N}^l
\gamma_{ik}\gamma_{jl}$. Also, \beq \label{eq:ineq} N^i N^jN^kN^l \gamma_{ik} \gamma_{jl}
= [N^iN^k\gamma_{ik}]^2 = [\delta^{ik}\gamma_{ik} - {Q}^{ik}\gamma_{ik}]^2
= [{Q}^{ik}\gamma_{ik}]^2 \leq {Q}^{ij} {Q}^{kl} \gamma_{ik} \gamma_{jl}, \eeq
since for a symmetric matrix $M$ we have $[\Tr(M)]^2 \leq c\Tr(M^2)$.
From (\ref{eq:ineq}), \bea {Q}^{ij} {Q}^{kl}\gamma_{ik}\gamma_{jl}
+ {Q}^{ij}{N}^k {N}^l \gamma_{ik}\gamma_{jl} &+& {N}^i {N}^j {Q}^{kl}  \gamma_{ik}
\gamma_{jl} +  {N}^i {N}^j {N}^k {N}^l\gamma_{ik} \gamma_{jl} \\ &
\leq & {Q}^{ij} {Q}^{kl}\gamma_{ik}\gamma_{jl} + 2{Q}^{ij} {N}^k {N}^l
\gamma_{ik}\gamma_{jl} + c{Q}^{ij} {Q}^{kl} \gamma_{ik} \gamma_{jl} \\
& \leq & 2c {Q}^{ij}({Q}^{kl}+{N}^k {N}^l) \gamma_{ik}\gamma_{jl} \\
&=& 2c{Q}^{ij}\delta^{kl} \gamma_{ik}\gamma_{jl}.\eea Using the fact
that $\gamma = \zeta  \de \alpha - \zeta \beta$ we have \beq
\label{eq:jan08081456} {Q}^{ij}\delta^{kl}\gamma_{ik}\gamma_{jl} =
{Q}^{ij}\delta^{kl}(\zeta \de \alpha)_{ik}(\zeta \de \alpha)_{jl} +
{Q}^{ij}\delta^{kl}(\zeta \de \alpha)_{ik}\zeta \beta_{jl} +
{Q}^{ij}\delta^{kl}\zeta \beta_{ik}(\zeta \de \alpha)_{jl} +
{Q}^{ij}\delta^{kl}\zeta \beta_{ik}\zeta \beta_{jl}\eeq where the
second and third term can be controlled by $\ve |\zeta \na \alpha|^2
+ \frac{1}{\ve}|\zeta \div \alpha|^2$ and the fourth term can be
controlled by $|\zeta \div \alpha|^2$. The first term in
(\ref{eq:jan08081456}) can be controlled as follows:
\begin{eqnarray} {Q}^{ij}\delta^{kl} (\zeta \de \alpha)_{ik}(\zeta
\de \alpha)_{jl} & = & {Q}^{ij}\delta^{kl}(\zeta \pa_i \alpha_k +
\zeta \pa_k \alpha_i)(\zeta \pa_j \alpha_l + \zeta \pa_l \alpha_j)
\\ & = & \label{eq:jan08081500} {Q}^{ij}\delta^{kl}(\zeta \pa_i
\alpha_k)(\zeta \pa_j \alpha_l) + {Q}^{ij}\delta^{kl}(\zeta \pa_i
\alpha_k)(\zeta \pa_l \alpha_j) \\ & + & {Q}^{ij}\delta^{kl}(\zeta
\pa_k \alpha_i)(\zeta \pa_j \alpha_l) + {Q}^{ij}\delta^{kl}(\zeta
\pa_k \alpha_i)(\zeta \pa_l \alpha_j). \end{eqnarray} Let $\na_{Q}^i
= {Q}^{ij} \pa_j$. Since ${Q}^{ij} = \delta_{mn}{Q}^{im}{Q}^{jn}$,
the first term in (\ref{eq:jan08081500}) can be bounded by $|\zeta
\na_{Q}[\alpha]|^2$. The second and third term in
(\ref{eq:jan08081500}) can be bounded by $\ve|\zeta \na \alpha|^2 +
\frac{1}{\ve}|\zeta \na_{Q}[\alpha]|^2$. The fourth term we
manipulate as follows: ${Q}^{ij}\delta^{kl}(\zeta \pa_k
\alpha_i)(\zeta \pa_l \alpha_j) = \delta_{mn} {Q}^{mi}(\zeta \pa_k
\alpha_i){Q}^{nj}(\zeta \pa^k \alpha_j)$ and \begin{eqnarray}
{Q}^{mi}(\zeta \pa_k \alpha_i) & = & {Q}^{mi}(\zeta \pa_i
\alpha_{k}) + {Q}^{mi}[\zeta \pa_{k} \alpha_i - \zeta \pa_i
\alpha_{k}] \\ & = & \label{eq:dec51136} \zeta
\na_{{Q}}^m[\alpha_k]  + {Q}^{mi}(\zeta \curl \alpha)_{{k}i}.
\end{eqnarray} Thus the fourth term in (\ref{eq:jan08081500}) can be
controlled by $(1 + \frac{1}{\ve})|\zeta \na_{Q}[\alpha]|^2 + |\zeta
\curl \alpha|^2 + \ve|\zeta \na \alpha|^2$. This concludes the
proof. \ep \bigskip

From lemma \ref{normalderivnoninvII} we have the following result:

\begin{lemma} \label{april04080805II} For $1 \leq s \leq 5$,
\beq \label{eq:l2} \|\zeta \alpha\|_{H^s(\Om_t)} \leq P\big[\|x\|_5\big]
\left[\|\zeta \alpha\|_{L^2(\Om_t)} +  \|\zeta \curl \alpha\|_{H^{s-1}(\Om_t)}
+ \|\zeta \div \alpha\|_{H^{s-1}(\Om_t)} + \sum_{j=1}^s\|\zeta \pt^j \alpha\|_{L^2(\Om_t)}\right]. \eeq
\end{lemma}

\bp The base case is when $s = 1$ we have on $\Om_t$, according to lemma \ref{normalderivnoninvII},
$\|\zeta \na \alpha\| \leq \|\zeta \curl \alpha\|_{L^2(\Om_t)} + \|\zeta \div \alpha\|_{L^2(\Om_t)}
+ \|\zeta \pt \alpha\|_{L^2(\Om_t)}$, which means that (\ref{eq:l2}) holds. Now suppose that
$s = 5$ and that we have the result for smaller $s$. Then, by lemma \ref{normalderivnoninvII}
we see that \begin{align} \label{eq:feb25081342} \|\zeta \na^5 \alpha\|_{L^2(\Om_t)}
& \leq \|\zeta \curl \alpha\|_{H^4(\Om_t)} + \|\zeta \div \alpha\|_{H^4(\Om_t)}
+ \|\zeta \pt \na^4 \alpha\|_{L^2(\Om_t)}. \end{align} To manipulate the second
to last term in (\ref{eq:feb25081342}) we write
\begin{align} \label{eq:feb25081343} \zeta \na^4 \pt \alpha - \zeta \pt \na^4 \alpha
= \sum (\na^j \pt x)(\zeta \na^{k + 1} \alpha) \end{align} summing over $j + k = 4$
such that $k \leq 3$. For $0 \leq j \leq 2$ we have
$\|\na^j\pt x\|_{L^\infty(\Om_t)} \leq \|x\|_5$ and we control
$\|\zeta \na^{k + 1} \alpha\|_{L^2(\Om_t)}$ by induction. For
$3 \leq j \leq 4$ we have $\|\na^j\pt x\|_{L^2(\Om_t)} \leq \|x\|_5$ and
$\|\zeta \na^2 \alpha\|_{L^\infty(\Om_t)} \leq P\big[\|x\|_5\|\big] \|\zeta \alpha\|_{H^4(\Om_t)}$
by Sobolev's inequality. This we control by induction. Now
\begin{align} \|\zeta \na^4\pt \alpha\|_{L^2(\Om_t)}
& \leq \|\zeta \curl \na^3\pt \alpha\|_{L^2(\Om_t)}
+ \|\zeta \div \na^3\pt \alpha\|_{L^2(\Om_t)} + \|\zeta \pt \na^3\pt \alpha\|_{L^2(\Om_t)} \\
& \label{eq:feb25081229} \leq \|\zeta \na^3[(\na \pt x)(\na \alpha)]\|_{L^2(\Om_t)}
+ \|\zeta \na^3\pt \curl \alpha\|_{L^2(\Om_t)} + \|\zeta \na^3\pt \div \alpha\|_{L^2(\Om_t)} \\
& + \|\zeta \pt \na^3\pt \alpha\|_{L^2(\Om_t)}. \end{align} The first term in (\ref{eq:feb25081229})
is controlled by
\begin{align} \label{eq:april04080831} \sum \|(\na^{j+1}\pt x)(\zeta \na^{k+1}
\alpha)\|_{L^2(\Om_t)} \end{align} where the sums is over $j+k = 3$. This term can be controlled
by $\|x\|_5\|\zeta \alpha\|_{H^4(\Om_t)}$. We control the second term in (\ref{eq:feb25081229}) by
\beq \label{eq:april04080848} \sum \|(\na^j \pt x)(\zeta \na^{k+1} \curl \alpha)\|_{L^2(\Om_t)} \eeq
where the sum is over $j+k=3$. We control this term by $\|x\|_5\|\zeta \curl \alpha\|_{H^4(\Om_t)}$.
Similarly for the third term in (\ref{eq:feb25081229}). \ep \bigskip

In this section we prove the following lemma.

\begin{lemma} \label{CSineqII} Let $\div \alpha$ and $\curl \alpha$ be defined as in lemma \ref{normalderivnoninvII}.
Then, for $s \leq 5$,
\begin{align} \label{eq:dec7810} \|\alpha\|_{H^s(\Om_t)} \leq P\big[\|x\|_5\big]\left[\|\alpha\|_{L^2(\Om_t)}
+ \|\div \alpha\|_{H^{s-1}({{\Om_t}})} + \|\curl \alpha\|_{H^{s-1}({{\Om_t}})}
+ \|(\ptt^{s-\frac12} \alpha) \cdot N\|_{L^2(\pOm_t)}\right], \end{align} where
$N$ is the outward unit normal to $\pOm_t$ and where $p(s)$ is a polynomial which depends on $s$.
Also, for $s \leq 5$, \begin{align} \label{eq:dec7811} \|\alpha\|_{H^s(\Om_t)}
\leq  P\big[\|x\|_5\big]\left[ \|\alpha\|_{L^2({{\Om_t}})} + \|\div \alpha\|_{H^{s-1}(\Om_t)}
+ \|\curl \alpha\|_{H^{s-1}(\Om_t)} + \|(\ptt^{s-\frac12} \alpha) \cdot Q\|_{L^2(\pOm_t)}\right] \end{align}
where $Q$ is a unit vector which is tangent to $\pOm_t$. \end{lemma}

\bp First we prove (\ref{eq:dec7810}) and (\ref{eq:dec7811}) for $s = 1$, then we will use lemma
\ref{april04080805II} to obtain the higher order results. Finally, we will use interpolation to
obtain the result for real $s$. Now \begin{align} \label{eq:april04081151} \|\na \alpha\|_{L^2(\Om_t)}
=  \int_{\Om_t} \pa_j \alpha_i \pa^j \alpha^i dx = \int_{\pOm_t}  \alpha_i N_j \pa^j \alpha^i dS(x)
- (\alpha, \Delta \alpha)_{\Om_t} \end{align} where we define $(\alpha, \Delta \alpha)_{\Om_t}
= \int_{\Om_t} \alpha_i \pa_j\pa^j \alpha^i dx$. And
\begin{align} - (\alpha, \Delta \alpha)_{\Om_t} &
= - \int_{\Om_t} \alpha_i \left[\pa^i \pa^j \alpha_j + \pa_j \pa^j \alpha^i
- \pa^i \pa^j \alpha_j\right]dx = \int_{\Om_t} \alpha_i \left[- \pa^i \div \alpha
+ \pa^j \left(\curl \alpha \right)^{i}_j\right]dx \\ & = - \int_{\pOm_t} N^i
\alpha_i \div \alpha dS(x) + \int_{\Om_t} \left[\div \alpha\right]^2 dx
+ \int_{\pOm_t} \alpha_i N^j \left(\curl \alpha \right)^{i}_jdS(x)
- \int_{\Om_t} \pa^j \alpha_i  \left(\curl \alpha \right)^{i}_jdx. \end{align}
Also, \begin{align} - \int_{\Om_t} \pa^j \alpha_i  \left(\curl \alpha \right)^{i}_jdx
= - \int_{\Om_t} \left(\curl \alpha \right)_i^j \left(\curl \alpha \right)^{i}_jdx
- \int_{\Om_t} \pa_i \alpha^j  \left(\curl \alpha \right)^{i}_jdx \end{align} and
\begin{align}- \int_{\Om_t} \pa_i \alpha^j  \left(\curl \alpha \right)^{i}_jdx =
- \int_{\pOm_t} N_i \alpha^j  \left(\curl \alpha \right)^{i}_jdS(x)
+ \int_{\Om_t} \alpha^j  \pa_i \left(\curl \alpha \right)^{i}_jdx. \end{align}
Moreover, \begin{align} \int_{\Om_t} \alpha^j  \pa_i \left(\curl \alpha \right)^{i}_jdx
& = \int_{\Om_t} \alpha^j \pa_i \left[\pa^i \alpha_j - \pa_j \alpha^i \right]dx \\ &
= (\alpha, \Delta \alpha)_{\Om_t} - \int_{\Om_t} \alpha^j \pa_i \pa_j \alpha^idx \\
& \label{eq:april04081148} = (\alpha, \Delta \alpha)_{\Om_t} - \int_{\pOm_t} \alpha^j N_j
\pa_i \alpha^i dS(x) + \int_{\Om_t} \left[\div \alpha\right]^2 dx. \end{align}
 From (\ref{eq:april04081148}) we see that
 \begin{align} - 2\left(\alpha, \Delta \alpha\right)_{\Om_t} & \label{eq:april04081150}
 = 2 \int_{\Om_t} \left[\div \alpha\right]^2 dx - \int_{\Om_t} \left(\curl
 \alpha \right)_i^j \left(\curl \alpha \right)^{i}_jdx \\ & - 2\int_{\pOm_t}
 \alpha \cdot N \div \alpha dS(x) + \int_{\pOm_t}\left( \alpha_i N^j -
 N_i \alpha^j\right)\left(\curl \alpha \right)^{i}_jdS(x). \end{align}
 The boundary terms from (\ref{eq:april04081151}) and (\ref{eq:april04081150}) are
 \beq \label{eq:bdryterms} \int_{\pOm_t}  \alpha_i N_j \pa^j \alpha^idS(x) - \int_{\pOm_t}
 \alpha \cdot N \div \alpha dS(x) + \frac12 \int_{\pOm_t}\left( \alpha_i N^j
 - N_i \alpha^j\right)\left(\curl \alpha \right)^{i}_jdS(x).\eeq The second term in
 (\ref{eq:bdryterms}) can be manipulated using $Q$: On $\pOm_t$,
 $\alpha = \alpha \cdot {N} {N} + Q\alpha$ and therefore
 \beq \label{eq:secondterm} - \alpha \cdot N \div \alpha = - (\pa_i N^i)
 [\alpha \cdot N]^2 - \alpha \cdot N \na_N [\alpha \cdot N] -
 \alpha \cdot N N_i\na_N[Q^i\alpha] - \alpha \cdot N \na_{Qi}[Q^i\alpha] \eeq
  where $\na_N = N^i \pa_i$. In the above, $- \alpha \cdot N N_i\na_N[Q^i\alpha]
= - [\alpha \cdot N]^2 N_i\na_N[Q^{ij}]N_j - \alpha \cdot N
N_i\na_N[Q^{ij}]Q_j\alpha$. And the third term from
(\ref{eq:bdryterms}) we manipulate as follows: \begin{eqnarray}
\frac12 \left( \alpha_i N^j - N_i \alpha^j\right)\left(\curl \alpha
\right)^{i}_j & = & \frac12\left( \alpha_i N^j - N_i
\alpha^j\left(\pa^i \alpha_j - \pa_j\alpha^i\right)\right) \\ & = &
\frac12\left[ \alpha_i N^j \pa^i \alpha_j -\alpha_i N^j
\pa_j\alpha^i - N_i \alpha^j\pa^i \alpha_j + N_i
\alpha^j\pa_j\alpha^i\right] \\ & = & \label{eq:thirdterm} \alpha_i
N^j \pa^i \alpha_j - \alpha_i N^j \pa_j\alpha^i. \end{eqnarray} The
second term in (\ref{eq:thirdterm}) cancels the first term in
(\ref{eq:bdryterms}). The first term in (\ref{eq:thirdterm}) we deal
with as follows: \beq \label{eq:thirdterm2} \alpha_i N^j \pa^i
\alpha_j = \alpha_i \pa^i[\alpha \cdot N] - \alpha_i \alpha_j (\pa^i
N^j) = \alpha\cdot N \na_N[\alpha \cdot N] + Q_i\alpha
\na_Q^i[\alpha \cdot N] - \alpha_i \alpha_j (\pa^i N^j).\eeq The
first term in (\ref{eq:thirdterm2}) cancels the second term in
(\ref{eq:secondterm}). The remaining terms therefore are
\begin{eqnarray}  \label{eq:dec1906} &&  \int_{{\pa {\Om_t}}}
\left[- (\pa_i N^i) [\alpha \cdot N]^2 - [\alpha \cdot N]^2
N_i\na_N[Q^{ij}]N_j - \alpha \cdot N
N_i\na_N[Q^{ij}]Q_j\alpha\right]dS(x) \\ & + & \int_{{\pa {\Om_t}}}
\left[- \alpha \cdot N \na_{Qi}[Q^i\alpha] + Q_i\alpha
\na_Q^i[\alpha \cdot N] - \alpha_i \alpha_j (\pa^i
N^j)\right]dS(x).\end{eqnarray} To get the lower order terms into
the form we want, we use the fact that we can trade normal and
tangential boundary components: Define $\tau_{ij} = 2
\alpha_i\alpha_j - \delta_{ij}(\alpha^k)(\alpha_k)$. Then \bea
\Big|\|\alpha \cdot N\|^2_{L^2(\pOm_t)} - \|\alpha \cdot
Q\|^2_{L^2(\pOm_t)} \Big| & = & \left|\int_{\pOm_t}
\left[N^iN^j\alpha_i\alpha_j -
Q^{ij}\alpha_i\alpha_j\right]dS(x)\right| \\ & = &
\left|\int_{\pOm_t} \left[2N^iN^j -
\delta^{ij}\right][\alpha_i\alpha_j]dS(x)\right| \\ & \leq & \sum_{k
= 1}^\mu \left|\int_{\pOm_t} \zeta_k N^iN^j\tau_{ij}dS(x)\right| \\
& \leq & \sum_{k = 1}^\mu \left|\int_{U_k \cap \Om_t} \zeta_k
(\pa^iN^j) \tau_{ij}dx\right| + \sum_{k = 1}^\mu \left|\int_{U_k
\cap \Om_t} \zeta_k N^j \pa^i \tau_{ij}dx\right|. \eea Now \bea
\pa^i \tau_{ij} & = & 2\div \alpha \alpha_j + 2 \alpha_i (\pa_i
\alpha_j) - (\pa_j \alpha_k)(\alpha^k) - (\alpha_k)(\pa_j \alpha^k)
\\ & = & 2\div \alpha \alpha_j + 2 \alpha_i (\pa_i \alpha_j) +
\left(- (\pa_j \alpha_k)(\alpha^k)  + (\pa_k
\alpha_j)(\alpha^k)\right) - (\pa_k \alpha_j)(\alpha^k) \\ &+&
\left(- (\alpha^k)(\pa_j \alpha^k)  + (\alpha^k)(\pa_k
\alpha^j)\right) - (\alpha^k)(\pa_k \alpha^j).\eea Thus, \beq
\label{eq:dec6642} \Big|\|\alpha \cdot N\|^2_{L^2({\pa {\Om_t}})} -
\|\alpha \cdot Q\|^2_{L^2(\pOm_t)} \Big| \leq \|x\|_5\left[\|
\alpha\|^2_{L^2(\Om_t)} + \|\div \alpha\|_{L^2(\Om_t)} + \|\curl
\alpha\|_{L^2(\Om_t)}\right].\eeq Hence all the lower order terms in
(\ref{eq:dec1906}) can be controlled by \[\|\alpha \cdot
N\|^2_{L^2(\pOm_t)} + \|\alpha\|^2_{L^2(\Om_t)}  + \|\div
\alpha\|^2_{L^2(\Om_t)} + \|\curl \alpha\|^2_{L^2(\Om_t)}\] or \beq
\label{eq:jan09081259} \|\alpha \cdot Q\|^2_{L^2(\pa {{\Om_t}})} +
\|\alpha\|^2_{L^2({{\Om_t}})}  + \|\div \alpha\|^2_{L^2({{\Om_t}})}
+ \|\curl \alpha\|^2_{L^2({{\Om_t}})}.\eeq To control the fourth
term in (\ref{eq:dec1906}) we use lemma \ref{algebra}: \begin{align}
\int_{\pOm_t} \alpha \cdot N \na_{Qi}[Q^i\alpha]dS(x) & =
\int_{\pOm_t} \alpha \cdot N c\pt[Q^i\alpha]dS(x) \\ & \leq
\left\|\alpha \cdot N \right\|_{H^\frac12(\pOm_t)} \|\alpha \cdot
Q\|_{H^\frac12(\pOm_t)}. \end{align} By using the fact that $ab \leq
\ve a^2 + \frac{b^2}{\ve}$ and the trace theorem we see that the
fourth term in (\ref{eq:dec1906}) can be controlled by both
\[\frac{1}{\ve}\|x\|_4^2\left[\left\|\alpha \cdot N
\right\|_{H^\frac12(\pOm_t)} + \| \alpha\|^2_{L^2({{\Om_t}})} +
\|\div \alpha\|^2_{L^2({{\Om_t}})} + \|\curl
\alpha\|^2_{L^2({{\Om_t}})}\right] + \ve
\|x\|_4^2\|\alpha\|^2_{H^1(\Om_t)}\] and
\[\frac{1}{\ve}\|x\|_4^2\left[\|\alpha \cdot Q\|_{H^\frac12(\pOm_t)}
+ \| \alpha\|^2_{L^2({{\Om_t}})} + \|\div
\alpha\|^2_{L^2({{\Om_t}})} + \|\curl
\alpha\|^2_{L^2({{\Om_t}})}\right]  + \ve
\|x\|_4^2\|\alpha\|^2_{H^1(\Om_t)}. \] The fifth term in
(\ref{eq:dec1906}) can be controlled similarly. This proves
(\ref{eq:dec7810}) and (\ref{eq:dec7811}) for $s=1$. We now prove
the estimate in terms of $\alpha \cdot N$, with the estimate in
terms of $\alpha \cdot Q$ following similarly.

Suppose that $2 \leq s \leq 3$ and that we have the result for
smaller $s$. Using lemma \ref{april04080805II}, we have
\begin{align} \|\alpha\|^2_{H^s(\Om_t)} & =
\|\alpha\|_{L^2(\Om_t)}^2 + \|\na \alpha\|_{H^{s-1}(\Om_t)}^2 \\ &
\leq P\big[\|x\|_5\big]\left[\|\alpha\|_{L^2(\Om_t)} + \|\na
\alpha\|_{L^2(\Om_t)} + \|\curl \alpha\|_{H^{s-1}(\Om_t)} + \|\div
\alpha\|_{H^{s-1}(\Om_t)} + \sum_{i = 1}^{s-1}\|\pa^i \na
\alpha\|_{L^2(\Om_t)}\right] \end{align} where $\pa$ denotes both
$\pt$ and $\eta \na$. Then \begin{align} \label{eq:april04081231}
\pa^{s-1} \na \alpha - \na \pa^{s-1} \alpha = \sum (\na^{i+1}
x)(\na^{j+1} \alpha) \end{align} where the sum is over $i+j=s-1$
such that $j \leq s-2$. The commutator in (\ref{eq:april04081231})
can be controlled by $\|x\|_5\|\alpha\|_{H^{s-1}(\Om_t)}$. Using the
above computation for the case $s=1$ we have \begin{align} \|\na
\pa^{s-1} \alpha\|_{L^2(\Om_t)} & \leq
P\big[\|x\|_5\big][\|\pa^{s-1} \alpha\|_{L^2(\Om_t)} + \|\div
\pa^{s-1} \alpha\|_{L^2(\Om_t)} + \|\curl \pa^{s-1}
\alpha\|_{L^2(\Om_t)} + \|(\pa^{s-1} \alpha) \cdot
N\|_{H^{\frac12}(\pOm_t)}] \\ & \leq
P\big[\|x\|_5\big][\|\alpha\|_{H^{s-1}(\Om_t)} + \|\div
\alpha\|_{H^{s-1}(\Om_t)} + \|\curl \alpha\|_{H^{s-1}(\Om_t)} +
\|(\pt^{s-1} \alpha) \cdot N\|_{H^{\frac12}(\pOm_t)}]. \end{align}
Now suppose that $s=4$. Then we control the commutator in
(\ref{eq:april04081231}) as follows: For $0 \leq i \leq 2$ we
estimate this term as above. For $i = 3$ we have $j=0$ and we can
control the commutator by $\|x\|_4\|\alpha\|_{H^3(\Om_t)}$. The case
for $s = 5$ follows similarly. Using corollary \ref{algebra2} we see
that \[\|(\pt^{s-1} \alpha) \cdot N\|_{H^{\frac12}(\pOm_t)} \leq
\|(\ptt^\frac12 \pt^{s-1} \alpha) \cdot
N\|_{L^2(\pOm_t)}\|N\|_{L^\infty(\pOm_t)} +
\|\alpha\|_{H^{s-1}(\pOm_t)}\|\ptt^2 N\|_{L^2(\pOm_t)}. \] Let
$\alpha$ have tangential Fourier expansion $\sum \alpha_k
e^{ik\theta}$. Then $\pt^{s-1} \alpha^j = \sum (ik)^{s-1} \alpha_k^j
e^{ik\theta}$ and $\ptt^\frac12 \pt^{s-1} \alpha^j = \sum \langle
k\rangle^\frac12 (ik)^{s-1} \alpha_k^j e^{ik\theta}$. Therefore
$(\ptt^\frac12 \pt^{s-1} \alpha) \cdot N = \sum \langle
k\rangle^\frac12 (ik)^{s-1} \alpha_k \cdot N e^{ik\theta}$. Thus
\begin{align} \|(\ptt^\frac12 \pt^{s-1} \alpha) \cdot
N\|_{L^2(\pOm_t)}^2 = \sum \langle k\rangle k^{2(s-1)} |\alpha_k
\cdot N|^2. \end{align} Also, $(\ptt^{s-\frac12} \alpha) \cdot N =
\sum \langle k\rangle^{s-\frac12} \alpha_k \cdot N e^{ik\theta}$.
Thus \begin{align} \|(\ptt^{s-\frac12} \alpha) \cdot
N\|_{L^2(\pOm_t)}^2 = \sum \langle k\rangle^{2(s-\frac12)} |\alpha_k
\cdot N|^2. \end{align} Thus $\|(\ptt^\frac12 \pt^{s-1} \alpha)
\cdot N\|_{L^2(\pOm_t)} \leq \|(\ptt^{s-\frac12} \alpha) \cdot
N\|_{L^2(\pOm_t)}$. By interpolation we now obtain the result for
non-integer $s$. This concludes the proof. \ep

\bibliographystyle{plain}

\end{document}